\documentclass{birkmult}

\newtheorem{Thm}{Theorem}
\newtheorem{Def}[Thm]{Definition}

\newtheorem{Rem}[Thm]{Remark}
\newtheorem{Prop}[Thm]{Proposition}

\newtheorem{Lem}[Thm]{Lemma}

\renewcommand{\theequation}{\thesection.\arabic{equation}}

\def \N{{\mathbb N}}
\def \C{{\mathbb C}}
\def \R{{\mathbb R}}

\def \Z{{\mathbb Z}}
\def\Box{\hfill\framebox(0.25,0.25){}}

\def \lam{\lambda}
\def\D{\Delta}
\def \cF{{\mathcal F}}
\def \al{\alpha}
\def \lam{\lambda}
\def \AA{{\mathcal A}}

\def\cS{{\mathcal S}}

\def\H{{\mathbb H}}

\def\supetage#1#2{
\sup_{\scriptstyle {#1}\atop\scriptstyle {#2}} }

\def\refer#1{~\ref{#1}}
\def\refeq#1{~(\ref{#1})}
\def\ccite#1{~\cite{#1}}

\def\eqdefa{\buildrel\hbox{\footnotesize def}\over =}

\let\wt=\widetilde

\begin{document}
\title[The heat kernel and frequency localized functions on the
Heisenberg group]{The heat kernel and frequency localized
functions \\ on the Heisenberg  group} 

\author[H. Bahouri]{Hajer Bahouri}
\address[H. Bahouri]%
{ Facult{\'e} des Sciences de Tunis\\ D{\'e}partement de
Math{\'e}matiques\\
 1060 Tunis\\TUNISIE }
\email{hajer.bahouri@fst.rnu.tn }
\author[I. Gallagher]{Isabelle Gallagher}
\address[I. Gallagher]%
{ Institut de Math{\'e}matiques de Jussieu UMR 7586\\
Universit{\'e} Paris VII\\ 175, rue du Chevaleret\\ 75013
Paris\\FRANCE } \email{Isabelle.Gallagher@math.jussieu.fr}

%\subjclass[2000]{Primary  ; Secondary }

\begin{abstract}
  The goal of this paper is to study the action of the heat operator on the Heisenberg group~$\H^d$, and in particular to characterize Besov spaces of negative index on~$\H^d$
 in terms of the heat kernel. That characterization   can be extended to positive indexes using Bernstein inequalities. As a corollary we obtain  a proof of refined Sobolev inequalities in~$\dot W^{s,p}$ spaces.
\end{abstract}

  \maketitle
  
%%%%%%%%%%%%%%%%%%%%%%%%%%%%%%%%%%%%%%%%%%%%%%%%%%%

\section{Introduction}
\setcounter{equation}{0} 
This paper is concerned mainly with a characterization of Besov spaces on the Heisenberg group using the heat kernel. In~\cite{bgx}, a Littlewood-Paley
decomposition on the Heisenberg group is constructed, and Besov
spaces are defined using that decomposition. It is classical that
in~$\R^{d}$  there is an equivalent definition, for negative
regularity indexes, in terms of the heat kernel. This
characterization  in~$\R^{d}$ can be extended to  positive regularity indexes
thanks to Bernstein's inequalities which express that 
derivatives act almost as homotheties on distributions, the Fourier
transform of which is supported in a ring of~$\R^d$ centered at zero.

The aim of this text is to present a similar characterization of
Besov spaces on~$\H^d$ using the heat flow. One of the main steps of the procedure in~$\R^d$ consists in studying frequency localized functions  and 
the action of derivatives, and more generally Fourier multipliers,  on such functions (the correponding inequalities for derivatives are known as Bernstein inequalities).  In the Heisenberg group there is a priori no simple notion of frequency localization, since the Fourier transform is a family of operators on a Hilbert space; however frequencies may be understood by studying the action of the Laplacian on a Hilbertian basis of that space, which allows to define a notion of frequency localization (see Definition~\ref{definlocfreqheis} below). One can then try to  investigate the
action of the semi-group of the heat equation on the Heisenberg
group on such frequency localized functions. That is achieved in this paper; we also  prove a similar characterization of Besov spaces in terms of the heat flow, as in the classical~$\R^d$ case. This allows to prove refined Sobolev inequalities, for~$\dot W^{s,p}$ spaces. Finally we are able by similar techniques to recover the  fact that the heat semi-group is the convolution by a   function in the Schwartz class (as in previous works by  Gaveau
in~\cite{gaveau} and Hulanicki in~\cite{hulanicki}).

Let us mention that by a different method, Furioli, Melzi and
Veneruso obtained in \cite{fmv} a characterization of Besov spaces
in terms of the heat kernel for Lie groups of polynomial growth.

\subsection{The Heisenberg group~$ \H^d$ }
In this   introductory section, let us recall some basic facts
on the  Heisenberg group~$ \H^d$. The Heisenberg group~$ \H^d
$ is the Lie group with underlying~$\C^d \times \R$ endowed with
the following product law:
\[ \forall \left((z, s), (z', s')\right)
\in {\H^d} \times {\H^d},  \quad (z, s)\cdot (z', s') =
(z+z', s+s'+ 2  \mbox{Im}(z \cdot \overline z')),  \]
where~$z\cdot \overline z'= \sum_{j=1}^d z_j  \overline z'_j$. It
follows that~$ \H^d$ is a non commutative group,  the identity of
which is~$(0,0)$;  the inverse of the element~$(z,s)$ is given
by~$(z,s)^{-1}= (-z,-s)$. The Lie algebra of left invariant vector
fields  on the Heisenberg group~$ \H^d $ is spanned by the
vector fields
\[Z_j =
\partial_{z_j} + i \overline z_j
\partial_{s}, \quad \overline Z_j = \partial_{\overline z_j} -  i
z_j
 \partial_{s}  \quad\hbox{and}\quad S= \partial_s=\frac {1 }{2i}[\overline
Z_j,Z_j],
\]
 with~$ j\in \{1,\dots,d\}$.
In all that follows, we shall denote by~${\mathcal Z}$ the family
of vector fields defined by~$ Z_j$ for~$j \in \{1,\dots, d\}$ and~$Z_{j }=\overline Z_{j-d}$
for~$j \in \{d+1,\dots, 2d\}$ and  for any multi-index~$\alpha \in
\{1,\dots, 2d\}^k $, we will write
\begin{equation}\label{prodcomp} {\mathcal Z}^{\alpha} \eqdefa
Z_{\alpha_1}\dots Z_{\alpha_k}.
\end{equation}
The space~$\H^{d}$ is endowed with a smooth  left invariant
measure, the Haar measure, which in the coordinate
system~$(x,y,s)$ is simply the Lebesgue measure~$dxdyds$.

Let us point out that on the Heisenberg group~$\H^d$, there is a
notion of dilation defined for~$ a > 0 $ by~$\delta_a ( z,  s ) =
( a z,  a^2s ). $
 The homogeneous dimension of~$\H^d$ is therefore~$N \eqdefa 2d +
 2$,  noticing that the Jacobian of the dilation~$\delta_a$ is~$a^N$.

The Schwartz space~${\mathcal S}(\H^d)$ on the Heisenberg group is defined as follows.
 \begin{Def}\label{defschwartz}
 The Schwartz space~${\mathcal S}(\H^d)$ is the set of smooth functions~$u$ on~$\H^d$ such that, for any~$k \in \N$, we have
 $$
 \|u\|_{k,{\mathcal S}} \eqdefa \supetage{|\alpha| \leq k}{(z,s) \in \H^d} \left| {\mathcal Z}^\alpha  \left((|
z|^2 -is))^{2k}
 u(z,s) \right)\right|
< \infty.$$
 \end{Def}
\begin{Rem}\label{defrho}
The  Schwartz space on the Heisenberg group~${\mathcal S}(\H^d)$ coincides with the classical Schwartz space~${\mathcal S}(\R^{2d+1})$. The weight in~$(z,s)$ appearing in the definition above is related to the fact that the Heisenberg distance to the origin is defined by~$\displaystyle \rho(z,s) \eqdefa \left(|
z|^4 + s^2\right)^\frac{1}4$.
\end{Rem}
Finally, let us present   the Laplacian-Kohn operator, 
    which is central in the study of partial differential equations on~$\H^d$, and is defined by 
     $$ \Delta_{ \H^d  }\eqdefa 2 \sum_{j=1}^{d}(Z_j \overline
Z_j  +  \overline Z_j Z_j ).
$$
Powers of that operator allow to construct positive order Sobolev spaces: for example 
we define the homogeneous space~$\dot W^{s,p}(\H^d)$, for $0<s < N/p$,  as the completion of~${\mathcal S}(\H^d)$ 
for the norm 
$$
\|f\|_{\dot W^{s,p}(\H^d)} \eqdefa \left\|
(-\Delta_{ \H^d  })^\frac{s}2 f
\right\|_{L^p(\H^d)}.
$$

\subsection{Statement of the results }
In~\cite{bgx} and~\cite{bg} a dyadic unity partition  is built on
the Heisenberg group~$\H^d $,  similar to the one defined in the classical~$\R^d$ case. A significant application of this decomposition
is the definition of Besov spaces on the Heisenberg group in the
same way as in the classical case (see~\cite{bgx},\cite{bg}). In Section~\ref{hn}, we shall
give a full account of this theory.

The main result of this paper  describes
the action of the semi-group associated with the heat equation on the
Heisenberg group, on a frequency localized function.  We
refer to Definition~\ref{definlocfreqheis} below for the notion of
a frequency localized function, which requires the definition of the Fourier transform on~$\H^d$, and is therefore slightly technical.

\begin{Lem}
\label{echantillonexpheis} {   Let~$(r_{1},r_{2})$ be two
positive real numbers, and define~${\mathcal C}_{(r_1,r_2)}=
{\mathcal C}(0,r_1,r_2) $  the ring centered at the origin, of
small and large radius respectively~$r_{1}$ and~$r_{2}$.  Two
positive constants~$c$ and~$C$ exist such that, for any real
number~$p \in [ 1,\infty]$, any couple~$(t,\beta)$ of positive
real numbers and  any function~$u$ frequency localized in the
ring~$\beta{\mathcal C}_{(\sqrt{r_1},\sqrt{r_2})}$, we have
\begin{equation} 
\label{freqlocheat} 
 \| e^{t \Delta_{\H^d}}u\|_{L^p(\H^d)}
 \leq Ce^{-ct\beta^2}
 \|u\|_{L^p(\H^d)}.
  \end{equation} }
\end{Lem}

That lemma is the  key argument in the proof of  the following
theorem which is well known in~$\R^{d}$ and proved by a different
method  in \cite{fmv}  for Lie groups of polynomial growth. The definition of Besov spaces is provided in the next section.
\begin{Thm}
\label{equivsobhomosmheis} {  Let $s$ be a  positive real number and
$(p,r)\in [1,\infty]^2$.  A  constant~$C$ exists which satisfies
the following property. For~$u \in \dot B^{-2s}_{p,r}(\H^d)$, we
have
\begin{equation} \label{equivnormheis}
C^{-1} \|u\|_{\dot B^{-2s}_{p,r}(\H^d)} \leq   \Bigl\|\|t^s
e^{t\Delta_{ \H^d  }} u\|_{L^p(\H^d)} \Bigr\|_{L^r(\R^+, \frac
{dt}{t})} \leq C\|u\|_{\dot B^{-2s}_{p,r}(\H^d)} .
\end{equation} }
\end{Thm}

\begin{Rem}
Thanks to Bernstein's inequalities (see Proposition~\ref{lech} below), we have $$\|u\|_{\dot B^{\sigma}_{p,r}(\H^d)}
\equiv \sup_{|\alpha|=k}\left\|(-\Delta_{\H^d})^\frac\alpha2 u\right\|_{\dot
B^{\sigma-k}_{p,r}(\H^d)}. $$ We deduce that the caracterization
of Besov spaces on the Heisenberg group in terms of the heat
kernel can be extended to any positive regularity index.
\end{Rem}
This characterization is useful for instance to prove refined Sobolev   inequalities. In this paper we will prove the following result.
\begin{Thm}
\label{refinedSobolev} 
Let~$p \in [1,\infty]$ and~$0<s<N/p$  be given. 
There exists a positive constant~$C$ such that for any  function~$f$ in~$\dot W^{s,p}(\H^d)$  we have
$$
\|f\|_{L^q(\H^d)} \leq C \|f\|_{\dot W^{s,p}(\H^d)}^{1-\frac{sp}N} \|f\|_{\dot B^{s-\frac{N}p}_{\infty,\infty}}^\frac{sp}N,
$$
with~$q = pN /(N-ps)$.
\end{Thm}
\begin{Rem}
This is a refined Sobolev inequality since it is easy to see that~$\dot W^{s,p}(\H^d) $ is continuously embedded in~$\dot B^{s-\frac{N}p}_{\infty,\infty}$, so that Theorem~\ref{refinedSobolev} is a refined version of the classical inequality
$$
\|f\|_{L^q(\H^d)} \leq C \|f\|_{\dot W^{s,p}(\H^d)}.
$$
The above continuous embedding is simply due to the following  estimate, applied to~$u = (-\Delta_{\H^d})^\frac{s}2 f$:
$$
\|u\|_{\dot B^{-\frac{N}p}_{\infty,\infty}}  = \sup_{t>0} t^{\frac N{2p}} \|e^{t \Delta_{\H^d}} u\|_{L^\infty(\H^d)}
\leq  C \|u\|_{L^p(\H^d)}.
$$
Note that in the special case when~$p = 2$, such an inequality was proved in~\cite{bg}, using the method developped in the classical case in~\cite{gmo}.
\end{Rem}

It turns out that the techniques involved in the proof of
Lemma\refer{echantillonexpheis} enable us to recover the following
theorem, which was proved (by different methods) by Gaveau
in~\cite{gaveau} and Hulanicki in~\cite{hulanicki}.
\begin{Thm}
\label{heatkernelheis}  {  There exists a function~$h \in
{\mathcal S}(\H^d)$ such that,  if~$u$ denotes  the solution of
the free heat equation on the Heisenberg group
\begin{equation}
\label{eq:ps} \left\{\begin{array}{rcl}
                \partial_t u- \Delta_{\H^d} u &= & 0
        \quad \mbox{in} \quad \R^+ \times \H^d, \\
     u_{|t=0}
    & = & u_0 ,
        \end{array}
        \right.
\end{equation}
then we have $$u(t,\cdot)=u_0\star h_t$$ where $\star$  denotes
the convolution on the Heisenberg group defined in
Section~\ref{hn} below, while~$h_{t}$ is defined by
 $$h_t(x,y,s)=\frac{1}{t^{d+1}}
h(\frac{x}{\sqrt{t}},\frac{y}{\sqrt{t}},\frac{s}{t}).$$}
\end{Thm}

The rest of this paper is devoted to the proof of
Theorems~\ref{equivsobhomosmheis} to~\ref{heatkernelheis}, as well as
Lemma~\ref{echantillonexpheis}.  

The structure of the paper is the following. First, in  Section~\ref{hn}, we
present  a short review of Littlewood-Paley theory on the
Heisenberg group, giving  the notation and results that will be
used in the proofs, as well as the main references of the theory.    Section~\ref{proof} is
devoted to the proof of Theorem~\ref{equivsobhomosmheis},
assuming Lemma~\ref{echantillonexpheis}, and finally the proof of
Lemma~\ref{echantillonexpheis} can be found in Section~\ref{lem}.
In Section~\ref{lem} we also  give the proofs of Theorems~\ref{refinedSobolev} and~\ref{heatkernelheis}.

%%%%%%%%%%%%%%%%%%%%%%%%%%%%%%%%%%%%%%%%%%%%%%%%%%%%%%%%%

\section{Elements of Littlewood-Paley theory on the Heisenberg group}
\setcounter{equation}{0} \label{hn}
\subsection{The Fourier transform on the Heisenberg group}
To introduce the Littlewood-Paley theory  on the Heisenberg group,
 we need to recall the definition of the Fourier transform in that
 framework. We refer for instance to~\cite{nach}, \cite{stein2} or~\cite{taylor}   for
 more details. The
Heisenberg group being non commutative,  the Fourier transform
on~$\H^d $ is defined using  irreducible unitary representations
of~$ \H^d $. 
As explained for instance
in~\cite{taylor} Chapter 2, all irreducible representations of~$
\H^d $ are unitarily equivalent to one of two
representations: the Bargmann representation or
the~$L^2$-representation.  The representations on~$L^2(\R^d)$ can
be deduced from Bargmann representations thanks to
 interlacing operators. We can consult J. Faraut and K. Harzallah\ccite {farautharzallah} for more details.
We shall  choose here the Bargmann representations
described by~$(u^\lam,  {\mathcal H}_{\lam} )$,  with~$\lam \in
\R \setminus \{0 \}$, where~${\mathcal H}_{\lam}$ are the spaces
defined by \[{\mathcal H}_{\lam} = \{F \mbox{ holomorphic on }
{\C}^d,  \| F\|_{{\mathcal H}_{\lam} }
 < \infty\},\]
 while we define
 \begin{equation}
 \label{normbargman}
\| F\|_{{\mathcal H}_{\lam} }^2 \eqdefa \left(\frac{2 |\lam |}{
\pi}
  \right)^d \int_{{\C}^d}
e^{-2|\lam| |\xi|^2 }|F(\xi)|^2 d\xi,
\end{equation}
 and~$u^\lam $ is the  map from~$ \H^d $ into
 the group of unitary operators of~${\mathcal H}_{\lam} $  defined by
 \[\begin{array}{c}
u^{\lam}_{z,  s } F(\xi) = F(\xi - \overline z) e^{i \lam s + 2
\lam (\xi \cdot  z - |z|^2/2)}  \quad \mbox{for} \quad \lam
>0, \\
u^{\lam}_{z,  s } F(\xi) = F(\xi - z) e^{i \lam s - 2 \lam (\xi
\cdot \overline z - |z|^2/2)}  \quad \mbox{for} \quad \lam < 0.
\end{array}
\]
 Let us notice that~${\mathcal H}_{\lam} $ equipped with the norm~(\ref{normbargman})
  is a Hilbert  space and that the monomials
$$
   F_{\al,  \lam} (\xi) = \frac{(\sqrt{2|\lam|} \: \xi)^\al}{\sqrt{\al !}} ,
\quad \al \in {\N}^d, $$ constitute an orthonormal basis.

If~$f $ belongs to~$ L^1(\H^d)$,  its Fourier transform is given by 
 \[ {\mathcal F}(f)(\lam) \eqdefa
\int_{\H^d} f(z, s) u^{\lam}_{z,  s } dz ds. 
\] Note that the
function~${\mathcal F}(f) $ takes its values in the bounded
operators on~$ {\mathcal H}_{\lam}$.  As in the~$ \R^d$ case, one
has a Plancherel Theorem and an inversion formula. More precisely, let~$ {\mathcal A}$  denote  the Hilbert
space of   one-parameter families~$ A = \{ A (\lam ) \}_{ \lam \in
\R \setminus \{0 \}}$  of operators on~$ {\mathcal H}_\lam $ which
are Hilbert-Schmidt for almost every~$\lam \in \R $ with  norm
\[ 
\| A \| = \left( \frac{2^{d-1}}{\pi^{d+1}} \int_{- \infty}^{\infty}
\| {A (\lam )}\|_{HS ({ {\mathcal H}_\lam })}^2 |\lam |^{d} d\lam
\right)^{\frac{1}{2}} <  \infty
\]
where~$\| {A (\lam )}\|_{HS ({ {\mathcal H}_\lam })} $ denotes the
Hilbert-Schmidt norm of the operator~$A(\lam )$.  Then the Fourier
transform can be extended to an isometry from~$ L^2( {\H}^d) $
onto~${\mathcal A }$ and we have the Plancherel formula:
$$ 
\| f \|_{ L^2( \H^d)
}^2  =
\frac{2^{d-1}}{\pi^{d + 1}} \sum_{\al \in \N^d} \int_{-
\infty}^{\infty} \| 
{\mathcal  F}(f)(\lam ) F_{\al,  \lam} \|_{
{\mathcal H}_\lam }^2  |\lam |^{d} d\lam.
$$
 On the other hand, if
\begin{equation}\label{injfouheis}
\sum_{\al \in \N^d} \int_{- \infty}^{\infty}
\| {{\mathcal F}(f)(\lam ) F_{\al, \lam} } \|_{{\mathcal
H}_{\lam}}  |\lam |^{d} d\lam < \infty
\end{equation}
then we have for almost every~$w$, 
\begin{equation}\label{inversionfourier}
f(w)=
\frac{2^{d-1}}{\pi^{d + 1}} \int_{- \infty}^{\infty} {\rm tr}
\left( u^{\lam}_{w^{-1}}{\mathcal F}(f)(\lam )\right)
  |\lam |^{d} d\lam,
  \end{equation}
  where
   $$ {\rm tr}
   \left(
   u^{\lam}_{w^{-1}}
    {\mathcal F}(f)(\lam )\right)=
   \sum_{\al \in \N^d} (u^{\lam}_{w^{-1}}{\mathcal F}(f)(\lam )F_{\al,  \lam}, F_{\al,
   \lam})_{{\mathcal H}_{\lam}}
   $$ denotes the trace of the operator~$u^{\lam}_{w^{-1}}{\mathcal F}(f)(\lam
   )$.
\begin{Rem}\label{remarkschwartz}The above hypothesis\refeq{injfouheis}
  is satisfied in~$ {\mathcal S} (\H^d)$, where~$ {\mathcal S} ({\H }^d)$ is defined in Definition~\ref{defschwartz}. This follows from Proposition~\ref{schwarzfourier} which is proved  for the sake of completeness, directly below its statement.
  \end{Rem}

Let us moreover
point out that  we have the following useful formulas, for any~$ k\in
 \{1,...,d\}$. 
 
 Denoting
 by~$1_{k} = (0,\dots,1,\dots)$ the vector whose~$k$ - component is one and all the others are zero,  one has  
  \begin{equation}\label{zk1}
  {\mathcal F} (Z_kf)(\lam)
F_{\al,  \lam} = - \sqrt{2|\lam|}\sqrt{\al_k +1}{\mathcal
F}(f)(\lam ) F_{\alpha + 1_{k}, \lam}
 \end{equation}
 if~$\lambda > 0$, and similarly
  \begin{equation}\label{zk2}
  {\mathcal F} (Z_kf)(\lam)
F_{\al,  \lam}= \sqrt{2|\lam|}\sqrt{\al_k }{\mathcal
F}(f)(\lam ) F_{\alpha - 1_{k}, \lam}
 \end{equation}
if~$\lambda < 0$. Furthermore, 
  \begin{equation}\label{zk3}
  {\mathcal F} (\overline{Z_k}f)(\lam)
F_{\al,  \lam}= \sqrt{2|\lam|}\sqrt{\al_k }{\mathcal
F}(f)(\lam ) F_{\alpha - 1_{k}, \lam}
 \end{equation}
 if~$\lambda > 0$, and 
  \begin{equation}\label{zk4}
  {\mathcal F} (\overline{Z_k}f)(\lam)
F_{\al,  \lam} =- \sqrt{2|\lam|}\sqrt{\al_k +1}{\mathcal
F}(f)(\lam ) F_{\alpha + 1_{k}, \lam}
 \end{equation}
 if~$\lambda < 0$.
 Therefore, we have easily, for any~$\rho \in \R$,
   \begin{equation}
\label{eq:deltah} {\mathcal F}((-\Delta_{ {\H}^d })^\rho f)(\lam)  F_{\al,
\lam} = \left(4|\lam| (2 |\al |+ d)\right)^\rho {\mathcal F}(f)(\lam )  F_{\al,
\lam}
\end{equation}
 and
\begin{equation*}\label{heatheiseq}
\cF \left(e^{t\Delta_{\H^d}} f \right)(\lam) F_{\al, \lam} =
e^{-t(4 |\lam| (2|\al| + d))} \cF (f)(\lam) F_{\al, \lam}.
\end{equation*}

Using those formulas, we can prove the following proposition, which justifies Remark~\ref{remarkschwartz} stated above. The proof of this proposition is new to our knowledge.
\begin{Prop}\label{schwarzfourier}
For any function~$f \in {\mathcal S}(\H^d)$, (\ref{injfouheis}) is satisfied. More precisely, for any~$\rho > \frac N 2$,
    there exists  a positive constant~$C$ such that
\[ \sum_{\al \in \N^d}\int_{-\infty}^{\infty}
\| 
{\mathcal F}(f)(\lambda )F_{\al,  \lambda}\|
_{{\mathcal H}_{\lambda}} | \lambda |^{d}\,
d \lambda \leq C \Bigl(\| f\|_{L^1({\H}^d)} + \|(- \Delta_{
{\H}^d })^\rho f\|_{L^1({\H}^d)}\Bigr).\]
\end{Prop}
Let us  prove that result. By definition of~$ {\mathcal S} (\H^d)$, for any~$\rho \in \R$, the
function~$(-\Delta_{ {\H}^d })^\rho f$ belongs to~$ {\mathcal S} (\H^d)$. Therefore, we can write, using~(\ref{eq:deltah}),
\[
 {\mathcal F}(f)(\lambda )F_{\al,  \lambda}
 = {\mathcal F}( (-\Delta_{ {\H}^d  })^{-\rho }  (-\Delta_{ {\H}^d })^ \rho f)(\lambda )F_{\al,  \lambda}  =(4| \lambda | (2 |\al |+ d))^{-\rho} {\mathcal F}(  (-\Delta_{ {\H}^d })^ \rho f)(\lambda )F_{\al,  \lambda}.
 \]
But that implies that
$$
\| {\mathcal F}(f)(\lambda )F_{\al,  \lambda}\|_{{\mathcal H}_{\lambda}}^2 =   (4| \lambda | (2 |\al |+
d))^{-2\rho}\Bigl(\frac{2 | \lambda |}{ \pi} \Bigr)^d   \int_{\C^d}
e^{-2| \lambda | |\xi|^2 }\Bigl |{\mathcal F}(  (-\Delta_{ {\H}^d })^\rho
f)(\lambda )F_{\al, \lambda}(\xi)\Bigr |^2 d\xi.
$$
According to the definition of the Fourier transform on the Heisenberg
group, we thus have
\begin{eqnarray*} 
\| {\mathcal F}(f)(\lambda )F_{\al,  \lambda}\|_{{\mathcal H}_{\lambda}}^2 
&=& (4|\lam| (2 |\al |+
d))^{-2\rho}\Bigl(\frac{2 |\lam |}{ \pi} \Bigr)^d  \int_{\C^d}
e^{-2|\lam| |\xi|^2 }  \\ & \times  &\!\! \biggl( \int_{{\H}^d} \!\!\!((-\Delta_{ {\H}^d })^\rho \! f(z,s))u_{z,s}^\lam F_{\al,\lam} dz ds
\overline{\int_{{\H}^d} \!\!\!((-\Delta_{ {\H}^d })^\rho\!
f(z',s'))u_{z',s'}^\lam F_{\al,\lam} dz' ds'\!\!}\biggr)\!d\xi.
\end{eqnarray*}
Fubini's theorem allows us to write
\begin{eqnarray*}\| {\mathcal F}(f)(\lambda )F_{\al,  \lambda}\|_{{\mathcal H}_{\lambda}}^2 
&=& (4|\lam| (2 |\al |+ d))^{-2\rho}
\\ & \times&\!\!\int_{{\H}^d}\!\!\int_{{\H}^d}\! (-\Delta_{ {\H}^d })^\rho f(z,s) 
  \overline{ (-\Delta_{ {\H}^d })^\rho
f(z',s')} (
 u_{z,s}^\lam F_{\al,  \lam}\!\mid\! u_{z',s'}^\lam F_{\al,  \lam})_{{\mathcal H}_{\lam}}
 dzdsdz'ds'.\end{eqnarray*}
Since the operators~$u_{z,s}^\lam$ and~$u_{z',s'}^\lam$ are
unitary on~${\mathcal H}_{\lam}$ and the family~$(F_{\al,  \lam})$ is a
Hilbert basis of~${\mathcal H}_{\lam}$, we deduce that $$\| {\mathcal F}(f)(\lambda )F_{\al,  \lambda}\|_{{\mathcal H}_{\lambda}} \leq  (4|\lam| (2 |\al |+
d))^{-\rho} \|(- \Delta_{ {\H}^d })^\rho f\|_{L^1({\H}^d)}.$$
To conclude we decompose the integral on~$\lambda$  into two parts, corresponding to ``high and low'' frequencies (the parameter~$|\lambda|^\frac12$ may be identified as a frequency, as will be clear in the next section -- it is in fact already apparent in~(\ref{eq:deltah}) above). Thus denoting~$\lam_m=(2m+d)\lam,$ we write
\begin{eqnarray*}& &  \sum_{\al \in \N^d} \int_{-\infty}^{\infty}
\| {\mathcal F}(f)(\lambda )F_{\al,  \lambda}\|_{{\mathcal H}_{\lambda}}  |\lam |^{d}
d\lam \leq  \sum_{m \in \N } \left(\begin{array}{c} m+d-1 \\ m
\end{array}\right) \Bigl( \|f\|_{L^1({\H}^d)} \int_{|\lam_m| \leq 1} |\lam|^{d}d\lam
 \\ & & +(4 (2 m+ d))^{-\rho}
\|(-\Delta_{ \H^d })^\rho f\|_{L^1({\H}^d)} \int_{|\lam_m| \geq
1}|\lam|^{-\rho}|\lam|^{d}d\lam \Bigr). \end{eqnarray*}
This gives the announced result
for $\rho>N/2.$  The proposition is proved. \Box

 Finally  the convolution product of two functions~$ f$ and~$g$ on~$ {\H^d}
$ is defined by
\[ f \star g ( w ) = \int_{\H^d} f ( w v^{-1} ) g( v) dv = \int_{\H^d} f ( v ) g(v^{-1} w)
dv.
\]
It should be emphasized that the convolution on the Heisenberg group  is not commutative.
Moreover if~$P$ is a left invariant vector field on~$\H^d$, then one sees easily that
\begin{equation}\label{Pfg}
P (f \star g) = f \star Pg,
\end{equation}
whereas in general~$P (f \star g) \neq Pf \star g$.
Nevertheless the usual Young inequalities are valid on the Heisenberg group, and one has moreover
\begin{equation}\label{fourconv}
{\mathcal F}( f \star g )( \lam ) = {\mathcal F}(f) ( \lam ) \circ
{\mathcal F}(g )( \lam ).
\end{equation}

\medskip

 It turns out that for radial functions on the Heisenberg group,
 the Fourier transform becomes simplified and puts into light the
 quantity that will play the role of the frequency size. Let us
 first recall the concept of radial functions on the Heisenberg
 group.
 \begin{Def}
{ A function~$f$ defined on the Heisenberg group~$\H^d$ is
said to be radial if it is invariant under the action of the
unitary group~$U(d)$ of~${\C}^d$, which means that for
any~$u \in U(d)$, we have
\[ f(z, s)= f(u(z), s), \quad \forall (z,s) \in \H^d.\]
 A radial function on the Heisenberg group can  then be written
under the form $$ f(z, s) = g(|z|, s).$$ }
\end{Def}
It can be shown (see for instance~\cite{nach}) that the
Fourier transform of radial functions of~$L^2(\H^d),$ satisfies
the following formulas: 
$$
{\mathcal F}(f)(\lam) F_{\al,  \lam} =
R_{|\al|} ( \lam ) F_{\al,  \lam}
$$
  where
   $$
   R_{m} (
\lam )= \left(\begin{array}{c} m+d-1 \\ m
\end{array}\right)^{-1} \int  e^{i\lam s}f(z, s) L_m^{(d-1)}(2|\lam|
|z|^2)e^{-|\lam| |z|^2}  dz ds,
  $$ 
  and where~$ L_m^{(p)}$ are
Laguerre polynomials  defined by
 $$
 L_m^{(p)} (t) = \sum_{k=0}^{m} ( -1 )^k \left(\begin{array}{c} m+p\\ m - k
 \end{array}\right ) \frac{t^k}{k!}, \quad t \geq  0, \quad m,  p \in \N.
 $$
Note that in that case
\begin{equation}\label{defL2NR}
\|f\|_{L^{2}(\H^{d})} = \|(R_{m})\|_{L^{2}_d(\N \times \R)} \eqdefa
\left(\frac{2^{d-1}}{\pi^{d+1}}\sum _m \left(\begin{array}{c}
m+d-1 \\ m
 \end{array}\right ) \int_{- \infty}^{\infty} |Q_m ( \lam )|^2 |\lam|^d
 d\lam  \right)^{\frac12},
\end{equation}
which corresponds to the Plancherel formula recalled above, in the radial case.  We also have the following inversion formula: if~$R_{m}$ belongs to~$L^{2}_d(\N \times \R)$ defined in~(\ref{defL2NR}), then the function
\begin{equation}\label{inversefourierradial}
f(z,s) = \frac{2^{d-1}}{\pi^{d+1}} \sum_m \int e^{-i \lam s}
R_m(\lam) L_m^{(d-1)} (2|\lam| |z|^2)e^{-|\lam| |z|^2}
|\lam|^d d\lam
\end{equation}
is a radial function in~$L^2(\H^d)$ and satisfies
$$
{\mathcal F}(f) (\lam) F_{\alpha, \lam} = R_{|\alpha|} (\lam) F_{\alpha, \lam} .
$$

\subsection{ Littlewood-Paley theory on the Heisenberg group}
Now we are ready to define the Littlewood-Paley decomposition on~$
\H^d$.
 We will not give any proof but refer to the construction
 in~\cite{bgx} and~\cite{bg} for all the details. We simply recall that
 the key point in the construction of the Littlewood-Paley decomposition on~$
\H^d$ lies in the following proposition proved in~\cite{bgx}. Note that Proposition~\ref{convergence} enables one to show in particular that functions of~$-\Delta_{\H^d}$ may be seen as convolution operators by Schwartz class functions (a result proved by Hulanicki~\cite{hulanicki} in the case of general nilpotent Lie groups).
\begin{Prop}
\label{convergence} {  For any~$ Q \in {\mathcal D}
 (\R \setminus \{0 \} )$, the series
$$
  g( z, s) = \frac{2^{d-1}}{\pi^{d+1}} \sum_m \int e^{-i \lam s}
Q((2m+d)\lam) L_m^{(d-1)} (2|\lam| |z|^2)e^{-|\lam| |z|^2}
|\lam|^d d\lam $$ converges in~${\mathcal S}  (\H^d)$.
  }
\end{Prop}

 The Littlewood-Paley operators are then constructed using the following proposition
  (see~\cite{bgx} and~\cite{bg}).
\begin{Prop}
\label{dyaheisenberg} { Define the ring~$ {\mathcal C}_0 = \left\{\tau \in \R, \: \frac{3}{4} \leq |\tau| \leq
\frac{8}{3}\right\}$   and   the ball~${\mathcal B}_0 = \left\{\tau \in \R, \: |\tau|
\leq \frac{4}{3}\right\} $. Then there exist two radial
functions~$\widetilde R^*$ and~$R^* $ the values of which are in
the interval~$[0,1]$, belonging respectively to~${\mathcal
D}({\mathcal B}_0)$ and to~${\mathcal D}({\mathcal C}_0)$  such
that $$ \forall \tau \in \R, \quad \widetilde R^*(\tau) + \sum_{j
\geq 0} R^*(2^{-2j}\tau) = 1 \quad  \mbox{and} \quad \forall  \tau \in \R^*,    \quad \sum_{j \in
\Z}
 R^*(2^{-2j}\tau)  =  1,
 $$
and satisfying as well the support properties $$|p-q|\geq 1
\Rightarrow \mbox{supp}~R^*(2^{-2q}\cdot)\cap \mbox{supp}~
R^*(2^{-2p}\cdot)=\emptyset$$ $$
 \mbox{and}\quad q\geq 1 \Rightarrow \mbox{supp}~
\widetilde R^*\cap \mbox{supp}~ R^*(2^{-2q}\cdot) = \emptyset. $$
Moreover, there  are radial functions of~$ {\mathcal
  S}(\H^d)$, denoted~$ \psi $ and~$\varphi$
such that $$ {\mathcal F}(\psi)(\lam)F_{\al, \lam} = \widetilde
R^*_{|\al|}(\lam)F_{\al, \lam}  \quad \mbox{and} \quad {\mathcal
F}(\varphi)(\lam)F_{\al, \lam} = R^*_{|\al|}(\lam)F_{\al, \lam},
$$ where we have noted~$\widetilde R^*_m(\tau) = \widetilde
R^*((2m+d) \tau)$ and~$R^*_m(\tau)  =
 R^*((2m+d)
\tau)$.
 }
\end{Prop}
 Now as in the~$ \R^d$ case, we define  
Littlewood-Paley operators in the following way.
\begin{Def}
The Littlewood-Paley operators~$ \Delta_j$ and~$ S_j$, for~$ j \in
\Z$,  are defined by
  \begin{eqnarray*} {\mathcal F}(\D_jf)(\lam)F_{\al, \lam}
&=&R^*_{|\al|}( 2^{-2j}\lam) {\mathcal F}(f)(\lam)F_{\al, \lam}
 \\ {\mathcal F}(S_jf)(\lam)F_{\al, \lam}& = &
\widetilde R^*_{|\al|}( 2^{-2j}\lam){\mathcal F}(f)(\lam)F_{\al,
\lam} . \end{eqnarray*}
\end{Def}
\begin{Rem} It is easy to see that $$ \Delta_ju=
u\star2^{Nj}\varphi(\delta_{2^{j}}\cdot)
 \quad \mbox{and} \quad
S_ju= u\star2^{Nj}\psi(\delta_{2^{j}}\cdot) $$ which implies that
 those operators map~$L^p$ into~$L^p$ for all~$ p \in
[1,\infty]$ with norms which do not depend on~$j$.
\end{Rem}
 Along the same lines as in  the~$ \R^d$ case,
 we can define homogeneous Besov spaces on
the Heisenberg group (see~\cite{bgx}).
\begin{Def}
\label{besovheis} {  Let~$ s \in \R$ be given, as well as~$ p$
and~$ r$, two real numbers in the interval~$[1,\infty] $. The
Besov space~$\dot B^{s}_{p,r}(\H^d)$ is the space of tempered
distributions~$u$ such that
\begin{itemize}
\item The series $ \sum^{m}_{-m} \Delta_{q}u$ converges
  to~$u$ in~${\mathcal S}'(\H^{d})$.
\item $ \|u\|_{\dot
B^{s}_{p,r}(\H^d)} \eqdefa \left\|2^{q s}\| \Delta_{q}
u\|_{L^{p}(\H^d)}\right\|_{\ell^r (\Z)}< \infty$.
\end{itemize}
}
\end{Def}
\begin{Rem}
Sobolev spaces~$\dot H^s(\H^d)$ have a characterization using Littlewood-Paley operators, as well as noninteger H\"older spaces  (see~\cite{bgx},\cite{bg}). More precisely, one has $\dot H^s(\H^d) = \dot
B^{s}_{2,2}(\H^d)$ for any~$s \in \R$, and for any~$\rho \in \R \setminus \N$,~$\dot C^\rho(\H^d) = \dot
B^{\rho}_{\infty,\infty}(\H^d)$.
\end{Rem}

 %%%%%%%%%%%%%%%%%%%%%%%%%%%%%%%%%%%%%%%%%%%%
 
\subsection{Frequency localized functions and Bernstein inequalities on the Heisenberg group}
 Let us first define the
concept of localization procedure in the frequency space in the
framework of the Heisenberg group. We will only state the
definition in the case of smooth functions -- otherwise one
proceeds by regularizing by convolution (see~\cite{bgx}
or~\cite{bg}).
\begin{Def}
\label{definlocfreqheis} {  Let~${\mathcal C}_{(r_1,r_2)}=
{\mathcal C}(0,r_1,r_2) $ be a ring  of~$\R$ centered at the
origin. A  function~$u$ in~$ {\mathcal S}(\H^d)$ is said to be
frequency localized in the ring~$2^{j}{\mathcal C}_{(\sqrt{ r_1},
\sqrt{r_2})}$, if
 \[ {\mathcal F}
(u)(\lam) F_{\al, \lam} = {\bf 1}_{(2 |\al| + d)^{-1}
2^{2j}{\mathcal C}_{(\sqrt{ r_1},
\sqrt{r_2})}}(\lam) {\mathcal F}(u)(\lam)F_{\al, \lam}.
\]
 }
\end{Def}
\begin{Rem}
Equivalently, a frequency localized function in the sense of Definition~\ref{definlocfreqheis} satisfies
$$
u = u \star \phi_j,
$$
where~$\phi_j =  2^{Nj}\phi(\delta_{2^{j}}\cdot)$, and~$\phi $ is a radial function in~${\mathcal S}(\H^d)$
such that
\[
{\mathcal F}(\phi)(\lam) F_{\al,  \lam} =
R((2 |\alpha|+d )\lam)  F_{\al,  \lam},
\]
with~$R$ compactly supported in a   ring of~$\R$ centered at zero.
\end{Rem}
In order to estimate  the cost of applying powers of the Laplacian on a frequency localized function, we shall need the following
proposition,  which ensures that the action powers of the Laplacian act as homotheties on such  frequency localized functions. The proof of that proposition may be found in~\cite{bg}.
 \begin{Prop}[\cite{bg}]
\label{lech}
{\sl
Let $p$ be an element of~$[1,\infty]$   
and let~$(r_{1},r_{2})$ be two
positive real numbers. Define~${\mathcal C}_{(r_1,r_2)}=
{\mathcal C}(0,r_1,r_2) $  the ring centered at the origin, of
small and large radius respectively~$r_{1}$ and~$r_{2}$. Then for any real number~$\rho$, there is a constant~$C_\rho$ such that if~$u $ is a function defined 
on~$\H^d$, 
 frequency localized in the
ring~$2^j{\mathcal C}_{(\sqrt{r_1},\sqrt{r_2})}$ then 
$$    C_\rho^{-1} 2^{-j \rho}  \| (-\Delta_{\H^d})^{\frac{\rho}{2}} u\|_{L^p(\H^d)} \leq \| u\|_{L^p(\H^d)} 
\leq C_\rho 2^{-j \rho}  \| (-\Delta_{\H^d })^{\frac{\rho}{2}} u\|_{L^p(\H^d)}.
$$}
\end{Prop}

  %%%%%%%%%%%%%%%%%%%%%%%%%%%%%%%%%%%%%%%%%%%%

\section{Proof  of Theorem~\ref{equivsobhomosmheis}}\label{proof}
\setcounter{equation}{0} In this section we shall prove
Theorem~\ref{equivsobhomosmheis}, assuming
Lemma~\ref{echantillonexpheis}. It turns out that the proof is
very similar to the~$\R^{d}$ case, and we sketch it here for the
convenience of the reader.

Let us start by estimating~$\|t^s e^{t\Delta_{\H^d}} u\|_{L^p}$.
Using Lemma~\ref{echantillonexpheis} and the fact that the
operator~$\Delta_j$ commutes with the operator~$e^{t\Delta_{{\bf
H}^d}}$, we can write
\[
\|t^s \Delta_j  e^{t\Delta_{\H^d}} u\|_{L^p}  \leq  Ct^s
2^{2js}e^{-ct2^{2j}} 2^{-2js}\|\Delta_j  u\|_{L^p} .
\]
Using   the definition of the homogeneous Besov  (semi) norm, we
get
\begin{eqnarray*}
\|t^s  e^{t\Delta_{\H^d}}  u\|_{L^p}
 & \leq  &  C \|u\|_{\dot B^{-2s}_{p,r}}  \sum_{j\in \Z} t^s
2^{2js}e^{-ct2^{2j}}c_{  r,j}
\end{eqnarray*}
where  ~$(c_{  r,j})_{j\in \Z}$ denotes, as in  all  this proof, a
generic element of the unit sphere of~$\ell^{  r}(\Z)$. In the
case when~$r=\infty$, the required inequality comes immediately
from  the following easy result: for any  positive~$s$, we have
\begin{equation}\label{caseinfinity}
\sup_{t>0} \sum_{j\in \Z} t^s  2^{2js}e^{-ct2^{2j}} <\infty .
\end{equation}
In the case when~$r<\infty$,  using the H\"older inequality with
the weight~$ 2^{2js}e^{-ct2^{2j}}$ and
Inequality~(\ref{caseinfinity}) we obtain
\begin{eqnarray*}
\int_0^\infty t^{rs} \|e^{t\Delta_{\H^d}}  u\|_{L^p}^r \frac {dt}
t & \leq  &   C \|u\|^r_{\dot B^{-2s}_{p,r}} \int_0^\infty \biggl(
\sum_{j\in \Z}  t^s 2^{2js}e^{-ct2^{2j}}\biggr)^{r-1} \biggl(
\sum_{j\in \Z}  t^s 2^{2js}e^{-ct2^{2j}}c_{r,j}^r\biggr) \frac
{dt} t \\ & \leq   & C \|u\|^r_{\dot B^{-2s}_{p,r}} \int_0^\infty
\sum_{j\in \Z}  t^s 2^{2js}e^{-ct2^{2j}}c_{r,j}^r \frac {dt} t
\cdotp
\end{eqnarray*}
This gives directly the result by Fubini's theorem.

In order to prove the other  inequality, let us observe that for
any~$s$ greater than~$-1$,  we have
\[
\int_0^\infty \tau^s  e^{-\tau} d\tau\eqdefa C_s.
\]
Using the  fact that the Fourier transform on the Heisenberg group  is injective,
we deduce the following identity (which may be easily proved by taking the Fourier transform of both sides)
\[
\Delta_j u = C_s^{-1} \int_0^\infty t^s(-\Delta_{\H^d})^{s+1}
e^{t\Delta_{\H^d }} \Delta_j u   dt.
\]
Then Lemma~\ref{echantillonexpheis}, the obvious
identity~$e^{t\Delta_{\H^d}}u= e^{\frac t 2 \Delta_{{\bf
H}^d}}e^{\frac t 2 \Delta_{\H^d}}u$ and the fact that the
operator~$\Delta_j$ commutes with the operator~$e^{t\Delta_{{\bf
H}^d}}$,  lead to
\begin{equation}\label{deltajulp}
\|\Delta_ju\|_{L^p}    \leq    C\int_0^\infty  t^s 2^{2j(s+1)}
e^{-ct2^{2j}} \| e^{t\Delta_{\H^d}}u\|_{L^p}  dt.
\end{equation}
In the case~$r=\infty$, we  simply write
\begin{eqnarray*}
\|\Delta_ju\|_{L^p}
 & \leq  &  C \Bigl(\sup_{t>0} t^s  \| e^{ t   \Delta_{\H^d}}u\|_{L^p}\Bigr)
\int_0^\infty  2^{2j(s+1)} e^{-ct2^{2j}} dt  \\ & \leq  &  C
2^{2js}\Bigl(\sup_{t>0} t^s   \| e^{ t \Delta_{{\bf
H}^d}}u\|_{L^p}\Bigr).
\end{eqnarray*}
In the case~$r<\infty$,   H\"older's inequality with the
weight~$e^{-ct2^{2j}}$ gives
\begin{eqnarray*}
\biggl(\int_0^\infty  t^s  e^{-ct2^{2j}} \| e^{t\Delta_{{\bf
H}^d}}u\|_{L^p} dt \biggr)^r & \leq   & \biggl(\int_0^\infty
e^{-ct2^{2j}}dt \biggr)^{r-1}\int_0^\infty  t^{rs} e^{-ct2^{2j}}
\| e^{t\Delta_{\H^d}}u\|_{L^p}^r dt\\ & \leq   &
C2^{-2j(r-1)}\int_0^\infty t^{rs} e^{-ct2^{2j}}  \|
e^{t\Delta_{\H^d}}u\|_{L^p}^r dt.
\end{eqnarray*}
Thanks to (\ref{caseinfinity}) and Fubini's theorem, we infer
from~(\ref{deltajulp}) that
\begin{eqnarray*}
\sum_j 2^{-2jsr}\|\Delta_ju\|_{L^p}^r   & \leq  &  C
 \int_0^\infty  \biggl(\sum_{j\in \Z} t2^{2j}   e^{-ct2^{2j}} \biggr)
t^{rs} \| e^{t\Delta_{\H^d}}u\|_{L^p}^r \frac{dt}  t\\ & \leq  & C
 \int_0^\infty  t^{rs} \|  e^{t\Delta_{\H^d}}u\|_{L^p}^r \frac{dt} t \cdotp
\end{eqnarray*}
The theorem is  proved. \Box

 %%%%%%%%%%%%%%%%%%%%%%%%%%%%%%%%%%%%%%%%%%%%
 
\section{Proofs of Lemma~\ref{echantillonexpheis} and Theorems~\ref{refinedSobolev} and~\ref{heatkernelheis} }
\setcounter{equation}{0}  \label{lem}  
Now we are left with the proof of Lemma~\ref{echantillonexpheis}, as well as Theorems~\ref{refinedSobolev}  and~\ref{heatkernelheis}.  
Lemma~\ref{echantillonexpheis} is proved in Paragraph~\ref{prooflemma}, 
while the proofs of Theorems~\ref{refinedSobolev}  and~\ref{heatkernelheis} can be found in 
  Paragraphs~\ref{proofrefined} and~\ref{proofthm3} respectively.

\subsection{ Proof of Lemma~\ref{echantillonexpheis} }\label{prooflemma} By density, it suffices to
suppose that the  function~$u$ is an element of~${\mathcal S}(\H^d)$. Now the frequency localization of~$u$ in the ring~$\beta
{\mathcal C}_{(\sqrt{r_1},\sqrt{r_2})}$ allows us to write
\begin{equation}
\label{estredheis}
 {\mathcal F} (e^{t\Delta_{\H^d}}u )(\lam) F_{\al,  \lam} =
  e^{-t\beta^{2}(4 | \beta^{-2}\lam | (2|\al| + d))} R_{|\al|}(\beta^{-2}\lam) {\mathcal F} (u)(\lam) F_{\al,
  \lam},
\end{equation}
with~$R_{|\al|}( \lam) = R((2 |\al| + d)\lam)$ and~$R \in {\mathcal
D}(\R \setminus \{0 \} )$ is equal to~$ 1 $ near the
ring~${\mathcal C}_{(r_1, r_2 )}$. We can then assume in what
follows that~$\beta=1$.

Since~$R $ belongs to~$ {\mathcal
D}(\R \setminus \{0 \} )$, Proposition~\ref{convergence} ensures the existence of a radial
function~$g^t \in {\cS}(\H^d)$ such that
\[
\cF (g^t)(\lam)F_{\al,  \lam}= e^{-t(4 | \lam| (2|\al| + d))}
R_{|\al|}(\lam)F_{\al,  \lam}.
\]
We deduce that$$e^{t\Delta_{\H^d}}u = u \star g^t .$$ If we prove
that two positive real numbers~$c$ and~$C$ exist such that, for
all positive~$t$, we  have
\begin{equation}
\label{echantillonexpdemo2} \|  g^t\|_{L^1({\H}^d)} \leq C
e^{-ct},
\end{equation}
then the  lemma is proved. To prove~(\ref{echantillonexpdemo2}), let us first recall that
thanks to Proposition~\ref{convergence}
\[ g^t( z, s) = \frac{2^{d-1}}{\pi^{d+1}} \sum_m \int e^{-i \lam s}
e^{-t(4 | \lam| (2m + d))} R((2 m + d)\lam) L_m^{(d-1)} (2|\lam|
|z|^2)e^{-|\lam| |z|^2} |\lam|^d d\lam.
\]
Now, we shall follow the idea of the proof of
Proposition~\ref{convergence} established in~\cite{bgx} to obtain
Estimate~(\ref{echantillonexpdemo2}). Let us denote by~${\mathcal
Q}$ the sub-space of~$ L^2_d ( \N \times \R ) $ (defined in~(\ref{defL2NR})) generated by the
sequences~$(Q_m)$ of the type 
\begin{equation}
\label{funspaceheis}Q_m(\lam)=\int_{\R^n}Q((2m+f(\sigma))\lam)P(\lam)
d\mu(\sigma) ,
\end{equation}
where~$\mu$ is a bounded measure compactly supported on~$\R^n$,
~$f$ is a bounded function on the support of~$\mu$, ~$P$ is a
polynomial function and~$ Q$ a function of~$  {\mathcal D}
 (\R \setminus \{0 \} )$ under the form
 \begin{equation}
\label{genspaceheis}Q(\tau)=  e^{-4t|\tau|} {\mathcal P}(t \tau) R(\tau),
\end{equation}
with~${\mathcal P}$ a  polynomial  and~$R $ a function of~$ {\mathcal D}
 (\R \setminus \{0 \} )$.

Now, let us recall the following useful formulas (proved for
 instance in \cite{bgx} and \cite{magnus}).
 \begin{Lem}
\label{deriv-Fourier}{    For any radial function~$ f \in \cS
(\H^d)$, we have for any~$m \geq 1$,
\begin{eqnarray*}
 \cF ((is-|z|^2)f)(m, \lam)& = & \frac{d}{d \lam}\cF {f}(m,  \lam )
 - \frac{m  }{\lam}\Bigl(\cF {f}(m,  \lam )- \cF {f}(m -1,  \lam
)\Bigr) \: 
 \mbox{for} \quad \lam >
 0  \quad \mbox{and}
 \\ \cF ((is-|z|^2)f)(m, \lam)& = & \frac{d}{d \lam}\cF {f}(m,  \lam ) +
 \frac{m + d }{|\lam |}\Bigl(\cF {f}(m,  \lam )- \cF {f}(m +1,  \lam )\Bigr) \: 
  \mbox{for} \quad \lam
 < 0.
\end{eqnarray*}
Moreover, we have the following classical property on Laguerre
polynomials : \begin{equation}\label{lagop}
|L^{(p)}_{m}(y)e^{-y/2} | \leq C_{p} (m +1 )^{p}, \quad \forall y
 \geq 0,
\end{equation}}
\end{Lem}

Let us start by proving that for any integer~$k$, one has the
following formula
\begin{equation} \label{formnoyheis} (is-|z|^2)^k g^t( z, s) =
\frac{2^{d-1}}{\pi^{d+1}} \sum_m \int e^{-i \lam s}
Q^{(k)}_{m}(\lam) L_m^{(d-1)} (2|\lam| |z|^2)e^{-|\lam| |z|^2}
|\lam|^d d\lam,\end{equation} where~$(Q^{(k)}_{m})$ is an element
of the space~${\mathcal Q}$. By induction the problem is reduced
to proving that for~$(Q_m)$ element of~${\mathcal Q}$, the
sequence~$(Q^{\star}_{m})$ defined as follows is still an element
of~${\mathcal Q}$: for all~$m \geq 1$
\begin{eqnarray*}
Q^{\star}_{m}(\lam) &= &\frac{d}{d\lam}Q_m ( \lam) -
\frac{m}{\lam}(Q_m(\lam) -Q_{m-1}(\lam)), \quad \lam >0,  \\
Q^{\star}_{m}(\lam)& = &\frac{d}{d\lam}Q_m ( \lam) + \frac{m +
d}{|\lam |}(Q_{m}(\lam) -Q_{m+1}(\lam)), \quad \lam < 0.
\end{eqnarray*}
Let us for instance compute~$Q^{\star}_{m}(\lam)$ for~$\lam >0$
and~$m\geq 1$. Considering \refeq{funspaceheis}, the Taylor
formula implies that
\[ \frac{m}{\lam}(Q_m(\lam) -Q_{m-1}(\lam))= 2m \int_{\R^n}\int_0^1
Q'((2m+f(\sigma)-2u)\lam)P(\lam)du d\mu(\sigma).\] Therefore
\begin{eqnarray*}
Q^{\star}_{m}(\lam) &= & \int_{\R^n} Q((2m+f(\sigma))\lam)P'(\lam)
d\mu(\sigma)\\ &+& \int_{\R^n} Q'((2m+f(\sigma))\lam)P(\lam)
f(\sigma)d\mu(\sigma)\\ &+& 2 \int_{\R^n}\int_0^1 \int_0^1
(2m+f(\sigma)-2us)\lam Q''((2m+f(\sigma)-2us)\lam)P(\lam)ududs
d\mu(\sigma)\\ &-& 2\int_{\R^n}\int_0^1 \int_0^1
Q''((2m+f(\sigma)-2us)\lam)\lam P(\lam)ududs
f(\sigma)d\mu(\sigma)\\ &+& 4 \int_{\R^n}\int_0^1 \int_0^1
Q''((2m+f(\sigma)-2us)\lam)\lam P(\lam)u^2dusds d\mu(\sigma).
\end{eqnarray*}
This proves that the sequence~$(Q^{\star}_{m})$ belongs to the
space~${\mathcal Q}$.

Now let us end the proof of Lemma~\ref{echantillonexpheis}:
defining
\[f_m^t(z,s) =  \int e^{-i \lam s}
Q_{m}(\lam) L_m^{(d-1)} (2|\lam| |z|^2)e^{-|\lam| |z|^2} |\lam|^d
d\lam, \] with~$(Q_{m})$ element of~${\mathcal Q}$, and in view of~(\ref{formnoyheis}) it is enough
to prove that  there exist two constants~$c$ and~$C$ which do  not
depend on~$m$, such that
\begin{equation}\label{conclusionheis}
|f^{t}_{m}(z,s)| \leq C
e^{-ct} \frac{1}{m^2}\cdotp
 \end{equation} Due to   the condition
on the support of the function~$R$ appearing
in\refeq{genspaceheis}, there exist  two fixed constants which
only depend on~$R$, denoted~$c_1$ and~$c_2$ such that 
$$
f_m^t(z,s)
= \int_{\R^n}\int_{c_1\leq |(2 m + f(\sigma))\lam)|\leq c_2} e^{-i
\lam s} 
 Q((2m+f(\sigma))\lam)P(\lam)
L_m^{(d-1)} (2|\lam| |z|^2)e^{-|\lam| |z|^2} d\mu(\sigma)
|\lam|^d d\lam. 
$$
 In view of\refeq{genspaceheis}
and\refeq{lagop}, we obtain 
$$
|f_m^t(z,s)| \leq
c_{d-1} \int_{\R^n}\int_{c_1\leq |(2 m + f(\sigma))\lam)|\leq c_2} e^{-ct}
m^{d-1} d\mu(\sigma) |\lam|^d d\lam, $$ which leads easily
to~(\ref{conclusionheis}) and ends the proof of the lemma. \Box

 \subsection{Proof of Theorem~\ref{refinedSobolev} }\label{proofrefined}
 The proof of Theorem~\ref{refinedSobolev} presented here relies on the maximal function on the Heisenberg group; before starting the proof let us collect a few useful results on this function, starting with the definition of the maximal function (the interested reader can consult~\cite{stein2} for details and proofs). 
 \begin{Def}
  Let~$f$ be in~$ L^1_{loc}(\H^d)$. The maximal function of~$f$
is defined by 
$$ 
Mf(z, s) \eqdefa \sup_{R>0} \frac{1}{m(B((z, s),
R))} \int_{B((z, s), R)} |f(z', s')| dz' ds', 
$$ where~$m(B((z, s), R))$ denotes the  measure of the Heisenberg
ball~$B((z, s), R)$ of center~$(z,s)$ and radius~$R$. 
\end{Def}
The key propoerties we will use on the maximal function are collected   in the following proposition.
\begin{Prop}\label{propmax}
The maximal function satisfies the following properties.
\begin{enumerate}
\item  If~$f$ is a function in~$L^p(\H^d)$, with~$1<p\leq \infty$, then~$Mf $ belongs to~$L^p(\H^d)$ and we have
$$
\|Mf\|_{L^p(\H^d)} \leq A_p \|f\|_{L^p(\H^d)},
$$
where~$A_p$ is a constant which  depends only on~$p$ and~$d$.

\item Let~$\varphi$ be a function in~$ L^1(\H^d)$ and suppose that the function~$\displaystyle \psi(w) \eqdefa 
\sup_{\rho(w') \geq \rho(w) } \varphi(w')$ belongs to~$L^1(\H^d)$,
where~$\rho$ denotes the Heisenberg distance to the origin defined in Remark~\ref{defrho}. Then for any measurable function~$f$, we have
$$ \Big|(f \star \varphi)(w)\Big| \leq
\|\psi\|_{L^1(\H^d)} Mf(w).
$$
\end{enumerate}
\end{Prop}

Now we are ready to prove Theorem~\ref{refinedSobolev}. By density we can suppose that~$f$ belongs to~${\mathcal S} (\H^d)$. Let us write
$$
f=\int_0^\infty e^{t\Delta_{\H^d}}\Delta_{\H^d} f dt
$$ 
and decompose the integral in two parts:
$$
f=\int_0^A e^{t\Delta_{\H^d}}\Delta_{\H^d} f dt+\int_A^\infty
e^{t\Delta_{\H^d}}\Delta_{\H^d} fdt,
$$ where~$A$ is a constant to
be fixed later.

On the one hand, by Theorem\refer{equivsobhomosmheis}, we have
$$
\|e^{t\Delta_{\H^d}}\Delta_{\H^d} f\|_{L^\infty}\leq
\frac{C}{{t}^{1+\frac{1}{2}(\frac N p-s)}}\|f\|_{ \dot B^{s-\frac
Np}_{\infty,\infty}(\H^d)}.
$$ 
Therefore after integration we get
 $$
 \int_A^\infty
\|e^{t\Delta_{\H^d}}\Delta_{\H^d}  f\|_{L^\infty}\leq A^{\frac
12(s-\frac Np)}\|f\|_{\dot B^{s-\frac Np}_{\infty,\infty}({\bf
H}^d)}.
$$ 
On the other hand, denoting by~$g=(-\Delta_{{\bf
H}^d})^{\frac s 2}f$, we have 
$$
e^{t\Delta_{{\bf
H}^d}}\Delta_{\H^d} f=\frac{ 1}{(-t)^{1-\frac s 2}} e^{t\Delta_{{\bf
H}^d}}(-t\Delta_{\H^d})^{1-\frac s 2} g .
$$ 
It is well-known that the heat
kernel on the Heisenberg group satisfies the second assumption of   Proposition~\ref{propmax}
(the reader can consult~\cite{bcg},~\cite{fmv} or~\cite{gaveau}), so
we deduce that
$$
\Big|e^{t\Delta_{\H^d}}(-t\Delta_{\H^d})^{1-\frac s
2}g(x)\Big|\leq C_s M_{g}(x), 
$$ 
where~$M_g(x)$ denotes the maximal
function of the function~$g$. This leads to 
$$ 
\Big| \int_0^A
e^{t\Delta_{\H^d}}\Delta_{\H^d} f dt \Big| \leq C A^{\frac s 2}
M_{g}(x).
$$ 
In conclusion, we get
 $$
 \Big|\int_0^\infty
e^{t\Delta_{\H^d}}\Delta_{\H^d} f(x)dt\Big|\leq C \Bigl(A^{\frac s
2}M_g(x)+A^{\frac 12( s-\frac N p)}\|f\|_{\dot B^{s-\frac
Np}_{\infty,\infty} (\H^d)}\Bigr), 
$$ 
and the choice of~$A$ such that~$\displaystyle A^{\frac N{2p}}M_g(x) = \|f\|_{\dot B^{s-\frac Np}_{\infty,\infty}
}  $ ensures that
$$
\Bigl|\int_0^\infty e^{t\Delta_{\H^d}}\Delta_{\H^d} f(x)dt\Bigr|
\leq C M_{g}(x)^{1-\frac{ps}{N}}\|f\|_{\dot B^{s-\frac
Np}_{\infty,\infty} (\H^d)}^{\frac{ps}{N}}.
$$
Finally taking the~$L^q$ norm
with~$q=\frac{pN}{N-ps}$, ends the proof of Theorem~\ref{refinedSobolev} thanks to Proposition~\ref{propmax}. \Box

 \subsection{Proof of Theorem~\ref{heatkernelheis} }\label{proofthm3} The proof of
Theorem~\ref{heatkernelheis} is similar to the proof of
Lemma~\ref{echantillonexpheis} and relies on the following result.
\begin{Lem}
\label{kernelheisenberg2} {   The series
\begin{equation}\label{seiesheiseq}h( z, s) = \frac{2^{d-1}}{\pi^{d+1}} \sum_m \int e^{-i \lam s}
e^{- 4 | \lam| (2m + d )}  L_m^{(d-1)} (2|\lam| |z|^2)e^{-|\lam|
|z|^2} |\lam|^d d\lam.
\end{equation} converges in~${\mathcal S}(\H^d)$.  }
\end{Lem}
Notice that Lemma \ref{kernelheisenberg2} 
implies directly the theorem, as by a rescaling, it is easy to see that the heat
kernel on the Heisenberg group is given by
\[ h_t(x,y,s)=\frac{1}{t^{d+1}}
h(\frac{x}{\sqrt{t}},\frac{y}{\sqrt{t}},\frac{s}{t})\cdot  \] \Box

 \noindent {\bf Proof of Lemma\refer{kernelheisenberg2}}
 Due to the sub-ellipticity of~$-\Delta_{\H^d}$ (see for instance \cite{bgx}),
 it suffices to prove that for any
integers~$k$ and~$\ell$,
 $$ 
\left \|(-\Delta_{\H^d})^{\ell} ( |z|^2-is
)^{k} h \right\|_{L^2 (\H^d)} < \infty .$$ In order to do so, let
us introduce  the set~$\wt {\mathcal Q}$ of sequences~$(Q_m)$ of
the type
\begin{equation}
\label{funspaceheis2}Q_m(\lam)=\int_{\R^n}Q((2m+\theta(\sigma))\lam)P(\lam)
d\mu(\sigma) ,
\end{equation}
where~$\mu$ is a bounded measure compactly supported on~$\R^n$,~$\theta$ is a bounded function on the support of~$\mu$,~$P$ is a
polynomial function and~$ Q $ a function of~$ {\mathcal
C}^{\infty} (\R \setminus \{0 \} )$ under the form
 \begin{equation}
\label{genspaceheis2}Q(\tau)= e^{-4|\tau|}{\mathcal P}(\tau),
\end{equation}
where~${\mathcal P}$ is a polynomial function. As in the proof of
Lemma\refer{echantillonexpheis}  and thanks to
Formula\refeq{eq:deltah} and Lemma\refer{deriv-Fourier}, we obtain
\[ \Delta_{\H^d}^{\ell} (  |z|^2-is )^{k}h( z, s) = \frac{2^{d-1}}{\pi^{d+1}} \sum_m \int e^{-i \lam s}
 Q_m^{\ell,k}(\lam) L_m^{(d-1)} (2|\lam| |z|^2)e^{-|\lam|
|z|^2} |\lam|^d d\lam,\] with~$(Q_m^{\ell,k})$ an element of~$\wt {\mathcal
Q}$ which ends the proof of the lemma thanks to\refeq{defL2NR}. \Box

 %%%%%%%%%%%%%%%%%%%%%%%%%%%%%%%%%%%%%%%%%%%%

\end{document}